\documentclass[a4paper,10pt]{article}

\usepackage[utf8]{inputenc}
\usepackage{amsmath,amssymb}
\newtheorem{thm}{Theorem}[section]
\newtheorem{lem}[thm]{Lemma}
\newtheorem{pro}[thm]{Proposition}
\newtheorem{cor}[thm]{Corollary}
\usepackage[all]{xy}

\title{Approximate Hermitian-Yang-Mills structures and semistability for Higgs bundles. \\
       II: Higgs sheaves and admissible structures.}
\author{S. A. H. Cardona\footnote{Electronic address: andres.holguin@cimat.mx}\\CIMAT A.C. - Via Jalisco S/N - 36240, Gto. - M\'exico}

\begin{document}

\maketitle

\begin{abstract}
We study the basic properties of Higgs sheaves over compact K\"ahler manifolds and we establish some results concerning the notion of 
semistability; in particular, we show that any extension of semistable Higgs sheaves with equal slopes is semistable. Then, we use the  
flattening theorem to construct a regularization of any torsion-free Higgs sheaf and we show that it is in fact a Higgs bundle. 
Using this, we prove that any Hermitian metric on a regularization of a torsion-free Higgs sheaf induces an admissible structure on the 
Higgs sheaf. Finally, using admissible structures we prove some properties of semistable Higgs sheaves.   
\end{abstract}

\section{Introduction}

As it is well known, in dimensions greater or equal than two, the Mumford-Takemoto semistability of a coherent sheaf makes reference
to its coherent subsheaves and not only to its subbundles \cite{Kobayashi}, \cite{Siu}, \cite{Lubke}. This is also the case for Higgs sheaves 
\cite{Simpson 2}, \cite{Bruzzo-Granha}, and hence the notion of semistability in higher dimensions makes reference to Higgs subsheaves. In this article, the basic properties of Higgs sheaves are studied; some of them are simple extensions to the Higgs 
case of classical results on coherent sheaves, however they play an important role in the theory. In particular, in the first part we show that 
the kernel and the image of any morphism of Higgs sheaves are Higgs sheaves. On the other hand, since the torsion subsheaf of any Higgs sheaf can be seen as a Higgs 
sheaf in a natural way, we have that the definition of stability (resp. semistability) can be written in terms of subsheaves with torsion-free quotients. This
in fact justifies the definition introduced by Simpson in \cite{Simpson} and that have been used in \cite{Cardona}. \\

From the classical theory of coherent sheaves, we know that the direct sum of semistable sheaves with equal slope is again semistable, we prove 
here that the same holds for semistable Higgs sheaves. Furthermore, there is an analog classical result concerning the tensor product. Namely, 
the tensor product of semistable sheaves is again semistable. Kobayashi \cite{Kobayashi} obtained this result in the particular case of holomorphic 
bundles over projective spaces; latter on, Simpson \cite{Simpson} proved this for Higgs bundles, again when $X$ is projective. The result for 
Higgs sheaves over compact K\"ahler manifolds has been proved recently by Biswas and Schumacher \cite{Biswas-Schumacher} using an extension 
of the Hitchin-Kobayashi correspondence for torsion-free Higgs sheaves. \\

The notion of admissible metric on a torsion-free sheaf was introduced by Bando and Siu \cite{Bando-Siu}, they used these metrics to prove a 
Hitchin-Kobayashi correspondence for stable sheaves, which is in fact a natural extension to coherent sheaves of the main result 
of Uhlenbeck and Yau in \cite{Uhlenbeck-Yau}. Biswas and Schumacher in \cite{Biswas-Schumacher}, introduced the notion of admissible 
Hermitian-Yang-Mills structure on a torsion-free Higgs sheaf and generlalized this correspondence to Higgs sheaves. Here, we review the notion of 
admissible structure on a torsion-free Higgs sheaf and we prove, in particular, that tensor products and direct sums of admissible structures 
are again admissible structures. Then, using Hironaka's flattening theorem \cite{Hironaka} we construct a regularization of a torsion-free 
Higgs sheaf and we show that it is in fact a Higgs bundle. As a consequence of this we prove that any Hermitian metric on a regularization of 
a torsion-free Higgs sheaf induces an admissible structure on this sheaf. \\

In \cite{Cardona}, it was studied the Donaldson functional and the notions of semistability and approximate Hermitian-Yang-Mills for Higgs bundles 
over compact K\"ahler manifolds and it was proved that, if the Donaldson functional of a Higgs bundle is bounded 
from below, then such a bundle admits an approximate Hermitian-Yang-Mills structure. Hence, by a result in \cite{Bruzzo-Granha} we know it is also 
semistable. Finally, it was shown that, at least in the one-dimensional case, the concept of approximate Hermitian-Yang-Mills structure is the 
differential-geometric counterpart of the notion of semistability. It was done using a decomposition of the Donaldson functional for Higgs 
bundles and following the Donaldson's proof in \cite{Donaldson-1}, which was inspired in fact by the classical works on Riemann surfaces of 
Atiyah and Bott \cite{Atiyah-Bott}, Hitchin \cite{Hitchin} and Narasimhan and Seshadri \cite{Narasimhan-Seshadri}. The problem of 
the existence of approximate Hermitian-Yang-Mills structures on semistable Higgs bundles has been recently studied in 
\cite{Jiayu-Zhang} using the Donaldson heat flow techinque and some non-linear analyisis. \\

In the final part of this article we introduce the notion of approximate Hermitian-Yang-Mills structure on a 
torsion-free Higgs sheaf over a compact K\"ahler manifold and we prove that this notion behaves well with respect to tensor products and direct sums; which 
is indeed a counterpart of a result on semistable Higgs sheaves. Finally, using the Hitchin-Kobayashi correspondence for Higgs sheaves we show 
that a restriction of a polystable Higgs sheaf remains polystable when we restrict it to certain open sets. As a 
consequence of this, it follows that the tensor product of polystable Higgs sheaves is again polystable. This, together with a result of 
Simpson \cite{Simpson 3} that guarantees the existence of Harder-Narasimhan filtrations for Higgs sheaves, implies the semistability of any 
tensor product of semistable Higgs sheaves. As we said before, this result has been proved by Biswas and Schumacher in \cite{Biswas-Schumacher}; 
here we present a different proof of their result. \\

\noindent{\bf Acknowledgements} \\

\noindent The author would like to thank his thesis advisor, Prof. U. Bruzzo, for his constant support and encouragement and also 
The International School for Advanced Studies (SISSA) at Trieste, Italy. The main part of this work has been 
done when it was a Ph.D. student there. The author wants to thank also the Centro de Investigaci\'on en Matem\'aticas (CIMAT A.C.) at Guanajuato, 
M\'exico. Finally, the author would like to thank O. Iena, G. Dossena and C. Mar\'in for some comments. 
Some results presented in this article were part of the author's Ph.D. thesis in the Mathematical-Physics sector at SISSA.

\section{Higgs sheaves}

Let $X$ be an $n$-dimensional compact K\"ahler manifold with K\"ahler form $\omega$, and let $\Omega_{X}^{1}$ be the cotangent sheaf to $X$, 
i.e., it is the sheaf of holomorphic one-forms on $X$. As it was defined in \cite{Cardona}, a Higgs sheaf ${\mathfrak E}$ over $X$ is a coherent 
sheaf $E$ over $X$, together with a morphism $\phi : E\rightarrow E\otimes\Omega_{X}^{1}$ of ${\cal O}_{X}$-modules (that is usually called the 
Higgs field), such that the morphism $\phi\wedge\phi : E\rightarrow E\otimes\Omega_{X}^{2}$ vanishes. \\

Using local coordinates on $X$ we can write $\phi=\phi_{\alpha}dz^{\alpha}$, where the index take values ${\alpha}=1,...,n$ and each 
$\phi_{\alpha}$ is an endomorphism of $E$. Then the condition $\phi\wedge\phi=0$ is equivalent to $[\phi_{\alpha},\phi_{\beta}]=0$ for all 
$\alpha,\beta$. This condition, also called the integrability condition, implies that the sequence 
\begin{equation}
 \xymatrix{
E \ar[r]  &  E\otimes\Omega_{X}^{1} \ar[r]  &   E\otimes\Omega_{X}^{2} \ar[r]  &  \cdots  
} \nonumber 
\end{equation}
naturally induced by the Higgs field is a complex of coherent sheaves. A Higgs subsheaf ${\mathfrak F}$ of ${\mathfrak E}$ is a 
subsheaf $F$ of $E$ such that $\phi(F)\subset F\otimes\Omega_{X}^{1}$, so that the pair ${\mathfrak F}=(F,\phi|_{F})$ becomes itself a 
Higgs sheaf. A Higgs sheaf ${\mathfrak E}=(E,\phi)$ is said to be torsion-free (resp. locally free, reflexive, normal, torsion) if the 
corresponding coherent sheaf $E$ is torsion-free (resp. locally free, reflexive, normal, torsion). A Higgs bundle ${\mathfrak E}$ is by 
definition a Higgs sheaf which is locally free, hence a Higgs line bundle is just a locally free Higgs sheaf of rank one.\\

Let ${\mathfrak E}_{1}$ and ${\mathfrak E}_{2}$ be two Higgs sheaves over a compact K\"ahler manifold $X$. A morphism 
between ${\mathfrak E}_{1}$ and ${\mathfrak E}_{2}$ is a map $f:E_{1}\longrightarrow E_{2}$ such that the diagram 
\begin{displaymath}
 \xymatrix{
E_{1} \ar[r]^{\phi_{1}} \ar[d]^{f}    &     E_{1}\otimes\Omega_{X}^{1} \ar[d]^{f\otimes 1}   \\
E_{2} \ar[r]^{\phi_{2}}               &     E_{2}\otimes\Omega_{X}^{1}   \\
}
\end{displaymath}
commutes. In the following we will write any morphism of Higgs sheaves simply as $f:{\mathfrak E}_{1}\longrightarrow{\mathfrak E}_{2}$. Now, 
let ${\mathfrak F}=(F,\phi_{F})$ be a Higgs subsheaf of ${\mathfrak E}=(E,\phi)$ and let $G=E/F$. Then, in particular   
\begin{displaymath}
 \xymatrix{
 0  \ar[r]   &   F  \ar[r]  &   E  \ar[r]  &  G \ar[r]  &  0 \\
} 
\end{displaymath}
is an exact sequence of coherent sheaves. Tensoring this by $\Omega^{1}_{X}$ we get the following exact sequence 
\begin{equation}
 \xymatrix{
F\otimes\Omega^{1}_{X} \ar[r]^{f}  &   E\otimes\Omega^{1}_{X} \ar[r]  &   G\otimes\Omega^{1}_{X} \ar[r] &  0\,. 
}   \nonumber
\end{equation}
Since $\Omega^{1}_{X}$ is locally free, the morphism $f$ is injective (see \cite{Griffiths-Harris}, Ch.V, for details) and one has the 
following commutative diagram    
\begin{displaymath}
 \xymatrix{
0 \ar[r]   &   F \ar[r] \ar[d]^{\phi_{F}}         &     E \ar[r] \ar[d]^{\phi}       &       G \ar[d]^{\psi} \ar[r]       &    0 \\
0  \ar[r]  &   F\otimes\Omega^{1}_{X} \ar[r]      &     E\otimes\Omega^{1}_{X} \ar[r]    &       G \otimes\Omega^{1}_{X} \ar[r]  &    0 \\
}
\end{displaymath}
in which the rows are exact. The morphism $\psi$ in the above diagram is defined by demanding that all diagram becomes commutative (it is in 
fact well-defined because the rows are exact). It follows from this that $\psi$ is a Higgs field for the quotient sheaf $G$ and we say that 
the Higgs sheaf ${\mathfrak G}=(G,\psi)$ is a Higgs quotient of $\mathfrak E$.\\  

The kernel and the image of morphisms of Higgs sheaves are Higgs sheaves. In fact, if $f:{\mathfrak E}_{1}\longrightarrow{\mathfrak E}_{2}$ 
is a morphism of Higgs sheaves, $K={\rm ker\,}f$ and $\iota:K\longrightarrow E_{1}$ denotes the obvious inclusion, we have the 
following commutative diagram (with exact rows)
\begin{displaymath}
 \xymatrix{
K   \ar[r]^{\iota} \ar[d]^{\phi}        &     E_{1} \ar[r]^{f} \ar[d]^{\phi_{1}}      &     E_{2} \ar[d]^{\phi_{2}}  \\
K\otimes\Omega^{1}_{X} \ar[r]^{\iota'}  &     E_{1}\otimes\Omega^{1}_{X} \ar[r]^{f'}  &     E_{2}\otimes\Omega^{1}_{X} \\
}
\end{displaymath}
where $\iota'=\iota\otimes 1$, $f'=f\otimes 1$ and $\phi$ is the restriction of $\phi_{1}$ to $K$. In that way the pair 
${\mathfrak K}=(K,\phi)$ becomes a Higgs subsheaf of ${\mathfrak E}_{1}$. Similarly, if $F={\rm im\,}f$, we denote by 
$j:F\longrightarrow E_{2}$ the inclusion morphism and write $f=j\circ p$, we obtain the following commutative diagram 
\begin{displaymath}
 \xymatrix{
E_{1} \ar[r]^{p}\ar[d]^{\phi_{1}}       &   F \ar[r]^{j} \ar[d]^{\psi}          &     E_{2} \ar[d]^{\phi_{2}}    \\
E_{1}\otimes\Omega^{1}_{X} \ar[r]^{p'}  &   F\otimes\Omega^{1}_{X} \ar[r]^{j'}  &     E_{2}\otimes\Omega^{1}_{X} \\
}
\end{displaymath}
where $p'=p\otimes 1$, $j'=j\otimes 1$ and $\psi$ is the restriction of $\phi_{2}$ to $F$. From this we get that 
${\mathfrak F}=(F,\psi)$ is a Higgs sheaf. Furthermore, from the above diagram it follows that $\mathfrak F$ is a Higgs 
subsheaf of ${\mathfrak E}_{2}$ and at the same time a Higgs quotient of ${\mathfrak E}_{1}$. \\  

A sequence of Higgs sheaves is a sequence of the corresponding coherent sheaves where each map is a morphism of Higgs sheaves. 
A short exact sequence of Higgs sheaves, also called an extension of Higgs sheaves or a Higgs extension 
\cite{Bruzzo-Granha}, \cite{Simpson 2}, is defined in the obvious way. \\ 

Let
\begin{equation}
 \xymatrix{
0 \ar[r]  &  {\mathfrak F} \ar[r]  &   {\mathfrak E} \ar[r]  &   {\mathfrak G} \ar[r] &  0  
}   \label{model exact seq sheaves}
\end{equation}
be an exact sequence of Higgs sheaves. Since in the ordinary case we identify $F$ with a subsheaf of $E$, we can see the Higgs field 
of $\mathfrak F$ as a restriction of the Higgs field of $\mathfrak E$, in that way we identify $\mathfrak F$ with a Higgs subsheaf of 
$\mathfrak E$. \\  

Let ${\mathfrak E}=(E,\phi)$ be a Higgs sheaf. The morphism $\phi$ can be considered as a section of ${\rm End\,}E\otimes\Omega^{1}_{X}$ 
and hence we have a natural dual morphism $\phi^{*}:E^{*}\rightarrow E^{*}\otimes\Omega^{1}_{X}$ and the pair 
${\mathfrak E}^{*}=(E^{*},\phi^{*})$ is a Higgs sheaf. Furthermore, if $Y$ is another compact K\"ahler manifold and $f:Y\rightarrow X$ is a 
holomorphic map, the pair defined by $f^{*}\mathfrak E= (f^{*}E,f^{*}\phi)$ is also a Higgs sheaf. \\

On the other hand, if ${\mathfrak E}_{1}=(E_{1},\phi_{1})$ and ${\mathfrak E}_{2}=(E_{2},\phi_{2})$ are Higgs sheaves, the pair
\begin{equation}
{\mathfrak E}_{1}\otimes{\mathfrak E}_{2}=(E_{1}\otimes E_{2}\,,\phi_{1}\otimes I_{2} + I_{1}\otimes\phi_{2})\,  \label{tensor product}
\end{equation}
where $I_{1}$ and $I_{2}$ are the identity endomorphisms on $E_{1}$ and $E_{2}$ respectively, is a Higgs sheaf. Additionally, 
if ${\rm pr}_{i}:E_{1}\oplus E_{2}\rightarrow E_{i}$ (with $i=1,2$) denote the natural projections and we define 
${\rm pr}^{*}_{1}\phi_{1}$ and ${\rm pr}^{*}_{2}\phi_{2}$ by ${\rm pr}^{*}_{1}\phi_{1}(v_{1},v_{2})=(\phi_{1}v_{1},v_{2})$ and 
${\rm pr}^{*}_{2}\phi_{2}(v_{1},v_{2})=(v_{1},\phi_{2}v_{2})$ where $v_{1}$ and $v_{2}$ are sections of $E_{1}$ and $E_{2}$ respectively, the pair
\begin{equation}
{\mathfrak E}_{1}\oplus{\mathfrak E}_{2}=(E_{1}\oplus E_{2}\,,\,{\rm pr}^{*}_{1}\phi_{1} + {\rm pr}^{*}_{2}\phi_{2})\,  \label{direct sum}
\end{equation}
is a Higgs sheaf.\\

In \cite{Cardona} it was defined the degree and the rank of $\mathfrak E$, denoted by ${\rm deg\,}{\mathfrak E}$ and ${\rm rk\,}{\mathfrak E}$ 
respectively, as the degree and rank of the coherent sheaf $E$. Hence, if ${\rm det}\,E$ denotes the determinant bundle of the coherent 
sheaf $E$ we have
\begin{equation}
 {\rm deg\,}{\mathfrak E} = \int_{X}c_{1}({\rm det}\,E)\wedge\omega^{n-1}\,.
\end{equation}

If the rank is positive, we introduce the quotient $\mu({\mathfrak E})={\rm deg\,}{\mathfrak E}/{\rm rk\,}{\mathfrak E}$ which is 
called the slope of the Higgs sheaf. A Higgs sheaf $\mathfrak E$ is said to be $\omega$-stable (resp. $\omega$-semistable) if it is 
torsion-free and for any Higgs subsheaf ${\mathfrak F}$ with $0<{\rm rk\,}{\mathfrak F}<{\rm rk\,}{\mathfrak E}$ we have the inequality 
$\mu({\mathfrak F})<\mu(\mathfrak E)$ (resp. $\le$). We say that a Higgs sheaf is $\omega$-polystable if it decomposes into a direct sum of two 
or more $\omega$-stable Higgs sheaves all these with the same slope. Hence, a $\omega$-polystable Higgs sheaf is $\omega$-semistable but not 
$\omega$-stable. Notice that the notion of stability (resp. semistability) for Higgs sheaves 
makes reference only to Higgs subsheaves. Then, in principle a sheaf could be stable (resp. semistable) in the Higgs case, but not in 
the ordinary case. Finally, since the degree and the rank of any Higgs sheaf is the same degree and rank of the corresponding coherent sheaf, 
from \cite{Kobayashi}, Ch. V, Lemma 7.3, we know that any exact sequence of Higgs sheaves (\ref{model exact seq sheaves}) satisfies
\begin{equation}
{\rm rk\,}{\mathfrak F}\,(\mu({\mathfrak E}) - \mu({\mathfrak F})) + {\rm rk\,}{\mathfrak G}\,(\mu({\mathfrak E}) - \mu({\mathfrak G}))=0\,. \label{basic prop. sheaves}
\end{equation}

In a similar way to the classical case, from this equality it follows that the condition of stability (resp. semistability) can be written in 
terms of quotient Higgs sheaves instead of Higgs subsheaves. To be precise, as a direct consequence of (\ref{basic prop. sheaves}) we obtain  
\begin{pro}
Let ${\mathfrak E}$ be a torsion-free Higgs sheaf over a compact K\"ahler manifold $X$ with K\"ahler form $\omega$. Then, ${\mathfrak E}$ is 
$\omega$-stable (resp. $\omega$-semistable) if $\mu({\mathfrak E})<\mu({\mathfrak G})$ (resp. $\le$) for every quotient Higgs sheaf 
${\mathfrak G}$ with $0<{\rm rk\,}{\mathfrak G}<{\rm rk\,}{\mathfrak E}$. 
\end{pro}

Let ${\mathfrak E}=(E,\phi)$ be a Higgs sheaf over $X$ and let $T$ be the torsion subsheaf of $E$, since\footnote{Notice that $\phi$ is 
a morphism of sheaves, then $\phi(T)$ is contained in the torsion of $E\otimes\Omega^{1}_{X}$, which is exactly $T\otimes\Omega^{1}_{X}$ 
because $\Omega^{1}_{X}$ is locally free.} $\phi(T)\subset T\otimes\Omega^{1}_{X}$, the pair ${\mathfrak T}=(T,\phi|_{T})$ is a Higgs subsheaf 
of $\mathfrak E$; in other words, the torsion subsheaf of a Higgs sheaf is always a Higgs sheaf. Therefore ${\rm deg\,}{\mathfrak T}\ge 0$, which 
implies that in the definition of stability (resp. semistability) we do not have to consider all quotient Higgs sheaves. Namely, we get 
\begin{pro}\label{Fundamental property of ss, Higgs sheaves}
Let ${\mathfrak E}$ be a torsion-free Higgs sheaf over a compact K\"ahler manifold $X$ with K\"ahler form $\omega$. Then \\
{\bf (i)} ${\mathfrak E}$ is $\omega$-stable (resp. $\omega$-semistable) if and only if $\mu({\mathfrak F})<\mu({\mathfrak E})$ (resp. $\le$) for any Higgs subsheaf 
${\mathfrak F}$ with $0<{\rm rk\,}{\mathfrak F}<{\rm rk\,}{\mathfrak E}$ and such that the quotient ${\mathfrak E}/{\mathfrak F}$ is torsion-free. \\
{\bf (ii)} ${\mathfrak E}$ is $\omega$-stable (resp. $\omega$-semistable) if and only if $\mu({\mathfrak E})<\mu({\mathfrak G})$ (resp. $\le$) for any torsion-free 
quotient Higgs sheaf ${\mathfrak G}$ with $0<{\rm rk\,}{\mathfrak G}<{\rm rk\,}{\mathfrak E}$.
\end{pro}
\noindent {\it Proof:} (i) and (ii) are clear in one direction. For the converse, suppose the inequality between slopes in (i) (resp. in (ii)) holds 
for proper Higgs subsheaves with torsion-free quotient (resp. for torsion-free quotients Higgs sheaves) and let us consider an exact sequence of Higgs sheaves 
as in (\ref{model exact seq sheaves}). \\

Let ${\mathfrak E}=(E,\phi)$ and denote by $\psi$ the Higgs field of $\mathfrak G$, that is, ${\mathfrak G}=(G,\psi)$. Now, let $T$ be the 
torsion subsheaf of $G$. As we said before, the Higgs field satisfies $\psi(T)\subset T\otimes\Omega^{1}_{X}$ and consequently the pair 
${\mathfrak T}=(T,\psi|_{T})$ is a Higgs subsheaf of $\mathfrak G$; with torsion-free Higgs quotient, say ${\mathfrak G}_{1}$. Then if we 
define ${\mathfrak F}_{1}$ by the kernel of the Higgs morphism ${\mathfrak E}\longrightarrow{\mathfrak G}_{1}$, we have the following 
commutative diagram of Higgs sheaves
\begin{displaymath}
 \xymatrix{
         &                                           &                               &      0  \ar[d]                         &         \\
         &     0   \ar[d]                            &                               &    {\mathfrak T} \ar[d]                &         \\
0 \ar[r] &  {\mathfrak F} \ar[r] \ar[d]              &  {\mathfrak E} \ar[r] \ar[d]^{\rm Id}  &  {\mathfrak G} \ar[r] \ar[d]  &  0      \\
0 \ar[r] & {\mathfrak F}_{1} \ar[r] \ar[d]           &    {\mathfrak E} \ar[r]       &  {\mathfrak G}_{1} \ar[d] \ar[r]       &  0      \\
         & {\mathfrak F}_{1}/ {\mathfrak F} \ar[d]   &                               &               0                        &         \\            
         &              0                            &                               &                                        &         \\
}
\end{displaymath}
in which all rows and columns are exact. From this diagram we have that $\mathfrak F$ is a Higgs subsheaf of ${\mathfrak F}_{1}$ with 
${\mathfrak T}\cong{\mathfrak F}_{1}/{\mathfrak F}$. Since $\mathfrak T$ is a torsion Higgs sheaf, ${\rm deg\,}{\mathfrak T}\ge 0$ and we obtain
\begin{equation}
 {\deg\,}{\mathfrak G} = {\deg\,}{\mathfrak T} + {\deg\,}{\mathfrak G}_{1} \ge  {\deg\,}{\mathfrak G}_{1}\,,  \nonumber
\end{equation}
\begin{equation}
 {\deg\,}{\mathfrak F}_{1} = {\deg\,}{\mathfrak F} + {\deg\,}{\mathfrak T} \ge  {\deg\,}{\mathfrak F}\,.  \nonumber
\end{equation}
Now, because $\mathfrak T$ is torsion we have ${\rm rk\,}{\mathfrak G}={\rm rk\,}{\mathfrak G}_{1}$ and 
${\rm rk\,}{\mathfrak F}_{1}={\rm rk\,}{\mathfrak F}$ and hence finally we get 
\begin{equation}
\mu({\mathfrak F})\le \mu({\mathfrak F}_{1})\,, \quad\quad\quad \mu({\mathfrak G}_{1})\le \mu({\mathfrak G})\,.   \nonumber
\end{equation}
At this point, the converse directions in (i) and (ii) follows from the last two inequalities and the fact that 
${\mathfrak G}_{1}$ is torsion-free.   \;\;Q.E.D.  \\

Let ${\mathfrak E}=(E,\phi)$ be a torsion-free Higgs sheaf of rank $r$ over a compact K\"ahler manifold $X$; from a classical result 
(see Ch.V in \cite{Kobayashi} for more details) we know that ${\rm det\,}E\cong(\bigwedge^{r}E)^{**}$, and hence via this isomorphism $\phi$ defines a Higgs 
field $\eta$ on ${\rm det\,}E$. Therefore the pair ${\rm det\,}{\mathfrak E}=({\rm det\,}E,\eta)$ is a locally free Higgs sheaf, i.e., the 
determinant bundle of a torsion-free Higgs sheaf is a Higgs bundle. \\

On the other hand (see \cite{Kobayashi}, Ch.V, Proposition 6.12), we know that the determinant bundle of any torsion-free sheaf $E$ 
satisfies $({\rm det\,}E)^{*}={\rm det\,}E^{*}$. Consequently, if $\mathfrak E$ is torsion-free, 
$\mu(\mathfrak E)=-\mu({\mathfrak E}^{*})$ and we have the following proposition, which is a natural extension to Higgs sheaves 
of a classical result. 
\begin{pro}\label{Imte, properties dual sheaves}
Let $\mathfrak E$ be a torsion-free Higgs sheaf over a compact K\"ahler manifold $X$ with K\"ahler form $\omega$. Then\\
{\bf (i)} If ${\rm rk\,}{\mathfrak E}=1$, then $\mathfrak E$ is $\omega$-stable.\\
{\bf (ii)} Let $\mathfrak L$ be a Higgs line bundle over $X$. Then ${\mathfrak L}\otimes{\mathfrak E}$ is $\omega$-stable (resp. $\omega$-semistable) if and only if 
$\mathfrak E$ is $\omega$-stable (resp. $\omega$-semistable).\\
{\bf (iii)} $\mathfrak E$ is $\omega$-stable (resp. $\omega$-semistable) if and only if ${\mathfrak E}^{*}$ is $\omega$-stable (resp. $\omega$-semistable). 
\end{pro}
\noindent {\it Proof:} (i) is a direct consequence of the definition of stability and (ii) is equal to the classical case. We prove here (iii) 
in the case of stability (the proof for semistability is similar and is obtained by replacing $<$ by $\le$ in the inequalities between 
slopes). \\

Assume first ${\mathfrak E}^{*}$ is $\omega$-stable and consider the exact sequence of Higgs sheaves (\ref{model exact seq sheaves}) with 
$\mathfrak G$ torsion-free. Dualizing it, we obtain the exact sequence 
\begin{equation}
 \xymatrix{
0 \ar[r]  &  {\mathfrak G}^{*} \ar[r]  &   {\mathfrak E}^{*} \ar[r]  &   {\mathfrak F}^{*}\,.   
}  \nonumber
\end{equation}
Now, since ${\mathfrak E}$ and $\mathfrak G$ are both torsion-free, from the above sequence we get that 
\begin{equation}
\mu({\mathfrak E})=-\mu({\mathfrak E}^{*})<-\mu({\mathfrak G}^{*})=\mu(\mathfrak G)\,, \nonumber
\end{equation}
and hence, by Proposition \ref{Fundamental property of ss, Higgs sheaves} it follows that $\mathfrak E$ is $\omega$-stable.\\

Now, assume that $\mathfrak E$ is $\omega$-stable and consider an exact sequence of Higgs sheaves
\begin{equation}
 \xymatrix{
0 \ar[r]  &  {\mathfrak F}' \ar[r]  &   {\mathfrak E}^{*} \ar[r]  &   {\mathfrak G}' \ar[r] &  0  
}  \nonumber
\end{equation}
with ${\mathfrak G}'$ torsion-free. Dualizing this sequence, we obtain again an exact sequence of Higgs sheaves 
\begin{equation}
 \xymatrix{
0 \ar[r]  &  {\mathfrak G}'^{*} \ar[r]  &   {\mathfrak E}^{**} \ar[r]  &   {\mathfrak F}'^{*}\,.   
}  \nonumber
\end{equation}

Now, the natural injection $\sigma: \mathfrak E \longrightarrow {\mathfrak E}^{**}$ defines $\mathfrak E$ as a Higgs subsheaf 
${\mathfrak E}^{**}$. From this and defining the Higgs sheaves ${\mathfrak H}'={\mathfrak E}\cap{\mathfrak G}'^{*}$ and 
${\mathfrak H}''={\mathfrak E}/{\mathfrak H}'$ we have the following commutative diagram:
\begin{displaymath}
 \xymatrix{
         &     0             \ar[d]                           &              0    \ar[d]                                  &                          &      \\
0 \ar[r] &  {\mathfrak H}' \ar[r] \ar[d]                      &  {\mathfrak E} \ar[r] \ar[d]^{\sigma}                     &  {\mathfrak H}'' \ar[r]  &    0 \\
0 \ar[r] & {\mathfrak G}'^{*} \ar[r] \ar[d]                &    {\mathfrak E}^{**} \ar[r] \ar[d]                 &  {\mathfrak F}'^{*}   &      \\
0 \ar[r] & {\mathfrak G}'^{*}/ {\mathfrak H}' \ar[r]\ar[d] &    {\mathfrak E}^{**}/ {\mathfrak E} \ar[r] \ar[d]  &  {\mathfrak T}'' \ar[r]  &    0 \\            
         &              0                                     &                0                                          &                          &         \\
}
\end{displaymath}
where ${\mathfrak T}''$ is the quotient of the injective morphism ${\mathfrak G}'^{*}/ {\mathfrak H}'\longrightarrow{\mathfrak E}^{**}/ {\mathfrak E}$; so that 
the sequence on the bottom becomes an exact sequence. \\

In the above diagram, all columns and arrows are exact and since $\mathfrak E$ is torsion-free, the quotient 
${\mathfrak E}^{**}/ {\mathfrak E}$ is a torsion sheaf supported on a set of codimension at least two, and hence, the same holds also for 
${\mathfrak G}'^{*}/ {\mathfrak H}'$ and ${\mathfrak T}''$. Therefore ${\rm deg\,}{\mathfrak G}'^{*}={\rm deg\,}{\mathfrak H}'$ and 
${\rm rk\,}{\mathfrak G}'^{*}={\rm rk\,}{\mathfrak H}'$. Consequently, ${\mathfrak G}'^{*}$ and ${\mathfrak H}'$ have the same slope and it 
follows
\begin{equation}
 \mu({\mathfrak G}')=-\mu({\mathfrak G}'^{*})=-\mu({\mathfrak H}')>-\mu(\mathfrak E)=\mu({\mathfrak E}^{*})\,,  \nonumber
\end{equation}
which means that ${\mathfrak E}^{*}$ is $\omega$-stable. \;\;Q.E.D.  \\ 
\begin{cor}\label{}
Let ${\mathfrak E}$ be a torsion-free Higgs sheaf over a compact K\"ahler manifold $X$ with K\"ahler form $\omega$. Then ${\mathfrak E}$ is $\omega$-stable 
(resp. $\omega$-semistable) if and only if the sheaf ${\mathfrak E}^{**}$ is $\omega$-stable (resp. $\omega$-semistable).
\end{cor}

The above Corollary is an immediate consequence of the part (iii) of Proposition \ref{Imte, properties dual sheaves}. It has been proved 
independently by Biswas and Schumacher (see \cite{Biswas-Schumacher}, Lemma 2.4 for details). 

\section{Semistable Higgs sheaves}

In a similar way to the classical case we have a simple result concerning the direct sum of semistable Higgs sheaves. Namely we have  
\begin{thm}\label{Direct sum ss sheaves}
Let ${\mathfrak E}_{1}$ and ${\mathfrak E}_{2}$ be two torsion-free Higgs sheaves over a compact K\"ahler manifold $X$ with K\"ahler form $\omega$. 
Then ${\mathfrak E}_{1}\oplus{\mathfrak E}_{2}$ is $\omega$-semistable if and only if ${\mathfrak E}_{1}$ and ${\mathfrak E}_{2}$ are both 
$\omega$-semistable with $\mu({\mathfrak E}_{1})=\mu({\mathfrak E}_{2})$.    
\end{thm}

\noindent {\it Proof:} Assume first that ${\mathfrak E}_{1}$ and ${\mathfrak E}_{2}$ are both $\omega$-semistable 
with $\mu({\mathfrak E}_{1})=\mu({\mathfrak E}_{2})=\mu$ and let $\mathfrak F$ be a Higgs subsheaf of ${\mathfrak E}_{1}\oplus{\mathfrak E}_{2}$. Then 
we have the following commutative diagram where the horizontal sequences are exact and the vertical arrows are injective
\begin{displaymath}
 \xymatrix{
0 \ar[r]   &   {\mathfrak F}_{1} \ar[r] \ar[d]    &     {\mathfrak F} \ar[r] \ar[d]                &     {\mathfrak F}_{2} \ar[d] \ar[r]  &    0 \\
0 \ar[r]   &   {\mathfrak E}_{1} \ar[r]   &     {\mathfrak E}_{1}\oplus{\mathfrak E}_{2} \ar[r]    &     {\mathfrak E}_{2} \ar[r]   &    0 \\
}
\end{displaymath}
where ${\mathfrak F}_{1}={\mathfrak F}\cap({\mathfrak E}_{1}\oplus 0)$ and ${\mathfrak F}_{2}$ is the image of $\mathfrak F$ under 
${\mathfrak E}_{1}\oplus{\mathfrak E}_{2}\longrightarrow{\mathfrak E}_{2}$. From the above diagram we have 
\begin{equation}
 {\rm deg\,}({\mathfrak E}_{1}\oplus{\mathfrak E}_{2}) = {\rm deg\,}{\mathfrak E}_{1} + {\rm deg\,}{\mathfrak E}_{2}\,.
\end{equation} 
Now, since by hypothesis ${\mathfrak E}_{1}$ and ${\mathfrak E}_{2}$ have the same slope $\mu$, we have 
$\mu({\mathfrak E}_{1}\oplus{\mathfrak E}_{2})=\mu$ and
\begin{equation}
 {\rm deg\,}{\mathfrak F}_{1}\le\mu\cdot{\rm rk\,}{\mathfrak F}_{1}\,, \quad\quad  {\rm deg\,}{\mathfrak F}_{2}\le\mu\cdot{\rm rk\,}{\mathfrak F}_{2}\,.  \nonumber
\end{equation}
From these inequalities we obtain 
\begin{equation}
 \mu(\mathfrak F) = \frac{{\rm deg\,}{\mathfrak F}}{{\rm rk\,}{\mathfrak F}} = \frac{{\rm deg\,}{\mathfrak F}_{1} + {\rm deg\,}{\mathfrak F}_{2}}{{\rm rk\,}{\mathfrak F}_{1} + {\rm rk\,}{\mathfrak F}_{2}}\le\mu \nonumber
\end{equation}
and the semistability of ${\mathfrak E}_{1}\oplus{\mathfrak E}_{2}$ is follows. \\

Conversely, suppose ${\mathfrak E}_{1}\oplus{\mathfrak E}_{2}$ is $\omega$-semistable. Since ${\mathfrak E}_{1}$ and ${\mathfrak E}_{2}$ are at 
the same time Higgs subsheaves and quotient Higgs sheaves of ${\mathfrak E}_{1}\oplus{\mathfrak E}_{2}$ we necessarily obtain
\begin{equation}
 \mu({\mathfrak E}_{1}\oplus{\mathfrak E}_{2}) = \mu({\mathfrak E}_{1}) = \mu({\mathfrak E}_{2})\,. \nonumber
\end{equation}
A Higgs subsheaf ${\mathfrak G}_{1}$ of ${\mathfrak E}_{1}$ is clearly a Higgs subsheaf of ${\mathfrak E}_{1}\oplus{\mathfrak E}_{2}$ and 
hence $\mu({\mathfrak G}_{1})\le \mu({\mathfrak E}_{1})$, which shows the semistability of ${\mathfrak E}_{1}$. A similar argument shows the 
semistability of ${\mathfrak E}_{2}$.  \;\;Q.E.D. \\

In the proof of the above result we showed also that the slope of ${\mathfrak E}_{1}\oplus{\mathfrak E}_{2}$ is the same slope of 
${\mathfrak E}_{1}$ and ${\mathfrak E}_{2}$, which says that the direct sum of semistable Higgs sheaves can never be stable. In fact,
even if they are stable, the direct sum is just polystable. \\

The definition of semistability for Higgs sheaves that we have introduced in the preceding section uses only proper Higgs subsheaves, this 
definition can be reformulated in terms of Higgs sheaves of arbitrary rank (non necessarily proper). Indeed, it was the way in which 
Kobayashi \cite{Kobayashi} introduced the notion of semistability for holomorphic vector bundles. Therefore, alternatively we can say that a torsion-free  
Higgs sheaf $\mathfrak E$ over a compact K\"ahler manifold $X$ is $\omega$-semistable if and only if $\mu(\mathfrak F)\le\mu(\mathfrak E)$ 
for every Higgs subsheaf $\mathfrak F$ with $0<{\rm rk\,}{\mathfrak F}\le{\rm rk\,}{\mathfrak E}$ or equivalently if and only if $\mu(\mathfrak E)\le\mu(\mathfrak Q)$ for 
every quotient Higgs subsheaf $\mathfrak Q$ with $0<{\rm rk\,}{\mathfrak Q}\le{\rm rk\,}{\mathfrak E}$. In fact, this equivalence is clear from the identity (\ref{basic prop. sheaves}). On the other hand, if we assume that a sheaf $\mathfrak E$ is semistable according to our original definition and $\mathfrak F$ 
is a Higgs subsheaf with ${\rm rk\,}{\mathfrak F}={\rm rk\,}{\mathfrak E}$, then we have a 
sequence 
\begin{equation}
 \xymatrix{
0 \ar[r]  &  {\mathfrak F} \ar[r]  &   {\mathfrak E} \ar[r]  &   {\mathfrak Q} \ar[r]  &    0  
}  \nonumber
\end{equation}
with $\mathfrak Q$ a torsion sheaf (it is a zero rank sheaf). From the above exact sequence it follows that 
${\rm deg\,}{\mathfrak E}= {\rm deg\,}{\mathfrak F} + {\rm deg\,}{\mathfrak Q}$ and since ${\rm deg\,}{\mathfrak Q}\ge0$, necessarily 
$\mu(\mathfrak F)\le\mu(\mathfrak E)$. This means that $\mathfrak E$ is semistable with respect to the new definition. The converse direction 
is immediate, so that the definition of semistability can be written in terms of Higgs subsheaves (or equivalently quotient Higgs sheaves) 
of arbitrary rank. \\

Using the above definition of semistability it is easy to 
prove the following, which is a natural extension to Higgs sheaves of a classical result of semistable sheaves. 

\begin{pro}\label{Property of morphisms ss}
 Let $f:{\mathfrak E}_{1}\longrightarrow{\mathfrak E}_{2}$ be a morphism of $\omega$-semistable (torsion-free) Higgs sheaves over a compact K\"ahler 
manifold $X$ with K\"ahler form $\omega$. Then we have the following:\\
{\bf (i)} If $\mu({\mathfrak E}_{1})>\mu({\mathfrak E}_{2}),$ then $f=0$ (i.e, it is the zero morphism).\\
{\bf (ii)} If $\mu({\mathfrak E}_{1})=\mu({\mathfrak E}_{2})$ and ${\mathfrak E}_{1}$ is $\omega$-stable, then 
           ${\rm rk}\,{\mathfrak E}_{1}={\rm rk}\,f({\mathfrak E}_{1})$ and $f$ is injective unless $f=0\,$.\\  
{\bf (iii)} If $\mu({\mathfrak E}_{1})=\mu({\mathfrak E}_{2})$ and ${\mathfrak E}_{2}$ is $\omega$-stable, then 
           ${\rm rk}\,{\mathfrak E}_{2}={\rm rk}\,f({\mathfrak E}_{1})$ and $f$ is generically surjective unless $f=0\,$.\\
\end{pro}
\noindent {\it Proof:} Assume that ${\mathfrak E}_{1}$ and ${\mathfrak E}_{2}$ are both $\omega$-semistable with slopes $\mu_{1}$ and $\mu_{2}$ 
and ranks $r_{1}$ and $r_{2}$ respectively, and let ${\mathfrak F}=f({\mathfrak E}_{1})$; then $\mathfrak F$ is a torsion-free quotient Higgs 
sheaf of ${\mathfrak E}_{1}$ and a Higgs subsheaf of ${\mathfrak E}_{2}$.\\

(i) Suppose that $\mu_{1}>\mu_{2}$ and $f\neq0$, then 
\begin{equation}
 \mu(\mathfrak F)\le\mu_{2}<\mu_{1}\le\mu(\mathfrak F)\,, \nonumber
\end{equation}
which is impossible. Therefore $f$ must be the zero morphism. \\

(ii) Assume $f\neq0$ and suppose that $\mu_{1}=\mu_{2}$ and ${\mathfrak E}_{1}$ is $\omega$-stable. If $r_{1}>{\rm rk\,}{\mathfrak F}$, then 
\begin{equation}
\mu(\mathfrak F)\le\mu_{2}=\mu_{1}<\mu(\mathfrak F)\,. \nonumber
\end{equation}
Hence, necessarily $r_{1}={\rm rk\,}{\mathfrak F}$ and $f$ is injective.\\
 
(iii) Assume $f\neq0$ and suppose that $\mu_{1}=\mu_{2}$ and ${\mathfrak E}_{2}$ is $\omega$-stable. If $r_{2}>{\rm rk\,}{\mathfrak F}$, then 
\begin{equation}
 \mu(\mathfrak F)<\mu_{2}=\mu_{1}\le\mu(\mathfrak F)\,, \nonumber
\end{equation}
and consequently $r_{2}={\rm rk\,}{\mathfrak F}$ and the result follows.   \;\;Q.E.D. \\

From Proposition \ref{Property of morphisms ss}, we have that any extension of semistable Higgs sheaves with the same slope must be 
semistable. Namely we have
\begin{cor}\label{extension of ss sheaves}
 Let 
\begin{equation}
 \xymatrix{
0 \ar[r]  &  {\mathfrak F} \ar[r]  &   {\mathfrak E} \ar[r]  &   {\mathfrak G} \ar[r] &  0  
}   \nonumber
\end{equation}
be an exact sequence of torsion-free Higgs sheaves over a compact K\"ahler manifold $X$ with K\"ahler form $\omega$. If $\mathfrak F$ 
and $\mathfrak G$ are both $\omega$-semistable and $\mu(\mathfrak F)=\mu(\mathfrak G)=\mu$, then $\mathfrak E$ is also $\omega$-semistable 
and $\mu(\mathfrak E)=\mu\,.$
\end{cor}

\noindent {\it Proof:} The fact that $\mu(\mathfrak E)=\mu$ follows from the identity (\ref{basic prop. sheaves}). Suppose that $\mathfrak E$ is not 
semistable and hence there exists a subsheaf $\mathfrak H$ destabilizing it, i.e., there exists a proper (non-trivial) Higgs subsheaf $\mathfrak H$ 
such that $\mu(\mathfrak H)>\mu$. Without loss of generality we can assume that $\mathfrak H$ is semistable\footnote{If it is not, 
we can destabilize $\mathfrak H$ with a Higgs subsheaf ${\mathfrak H}'$. If it is semistable we stop, if it is not, then we repeat this 
procedure. Clearly this finishes after a finite number of steps, since in the extreme case we get a Higgs sheaf of rank one, which is in 
particular stable by Proposition \ref{Imte, properties dual sheaves}, part (i).}. Then we have a morphism 
$f:{\mathfrak H}\longrightarrow{\mathfrak G}$ with $\mu({\mathfrak H})>\mu({\mathfrak G})$, and from Proposition \ref{Property of morphisms ss} we 
have $f=0$. Therefore, there exists a morphism $g:{\mathfrak H}\longrightarrow{\mathfrak F}$ where $\mu({\mathfrak H})>\mu({\mathfrak F})$, and we 
have again $g=0$, which means that ${\mathfrak H}$ must be trivial and from this we have a contradiction. \;\;Q.E.D. \\

\section{Admissible structures}

Bando and Siu \cite{Bando-Siu} introduced the notion of admissible metric on a coherent sheaf; they used these metrics to prove a Hitchin-Kobayashi 
correspondence for stable torsion-free sheaves. Admissible structures have been used in \cite{Biswas-Schumacher} and \cite{Cardona} to extend the results 
of Bando and Siu to Higgs sheaves and to construct a Donaldson's functional for these objects. Here we review the main definitions concerning to 
this notion.\\

Let us consider a torsion-free Higgs sheaf ${\mathfrak E}$ over a compact K\"ahler manifold $X$ and let $S=S(\mathfrak E)\subset X$ be its singularity set; 
as it is well known, $S$ is a complex analytic subset with ${\rm codim}(S)\ge2$. An {\it admissible structure} on $\mathfrak E$ is an Hermitian 
metric $h$ on the bundle ${\mathfrak E}|_{X\backslash S}$ such that: \\

\noindent{\bf (i)} The Chern curvature $R_{h}$ is square-integrable, and \\
{\bf (ii)} The mean curvature $K_{h}=i\Lambda R_{h}$ is $L^{1}$-bounded. \\   

Since $S({\mathfrak E}^{**})\subset S(\mathfrak E)$, an admissible structure on ${\mathfrak E}^{**}$ induces an admissible structure on 
$\mathfrak E$. Now, the converse is also true if we relax the notion of admissibility to make only reference to certain open sets. In fact, 
this notion can be modified in the following way: an admissible structure on a Higgs sheaf $\mathfrak E$ is an Hermitian metric $h$ defined on an open 
set $U$, such that $X\backslash U$ is a complex analytic subset of codimension at least two, which contains the singularity set of $\mathfrak E$. 
In this sense, any admissible structure on $\mathfrak E$ induces an admissible structure on ${\mathfrak E}^{**}$. From now on 
(if necessary) we will understand admissible structures in this modified version. 

\begin{pro}\label{tensor product of admissible h's}
 Let ${\mathfrak E}_{1}$ and ${\mathfrak E}_{2}$ be two torsion-free Higgs sheaves over a compact K\"ahler manifold $X$ and let $h_{1}$ and 
$h_{2}$ be two admissible structures on these Higgs sheaves. Then, $h_{1}\otimes h_{2}$ and $h_{1}\oplus h_{2}$ are admissible structures on 
the tensor product ${\mathfrak E}_{1}\otimes{\mathfrak E}_{2}$ and the Whitney sum ${\mathfrak E}_{1}\oplus{\mathfrak E}_{2}$, respectively. 
\end{pro}
\noindent {\it Proof:} Suppose $h_{1}$ and $h_{2}$ are admissible structures and let $S_{1}$ and $S_{2}$ the singularity sets of 
${\mathfrak E}_{1}$ and ${\mathfrak E}_{2}$, respectively. Then, $h_{1}$ and $h_{2}$ are Hermitian metrics on ${\mathfrak E}_{1}|_{U_{1}}$ and 
${\mathfrak E}_{2}|_{U_{2}}$ for some open sets $U_{1}$ and $U_{2}$ in $X$, where $X\backslash U_{1}$ and $X\backslash U_{2}$ are 
complex analytic subsets of codimension greater or equal than two, containing the sets $S_{1}$ and $S_{2}$ respectively. \\

Since $X\backslash U_{1}\cap U_{2}$ is the union of $X\backslash U_{1}$ and $X\backslash U_{2}$, it is a closed analytic subset of codimension at 
least two containing $S_{1}\cup S_{2}$. From this we obtain that $h_{1}\otimes h_{2}$ and $h_{1}\oplus h_{2}$ are Hermitian metrics on the Higgs 
bundles ${\mathfrak E}_{1}\otimes{\mathfrak E}_{2}|_{U_{1}\cap U_{2}}$ and ${\mathfrak E}_{1}\oplus{\mathfrak E}_{2}|_{U_{1}\cap U_{2}}$ 
respectively. \\
  
On the other hand, if $K_{1\otimes2}$ denotes the Chern mean curvature of $h_{1}\otimes h_{2}$, we have from classical identities 
(see \cite{Kobayashi}, Ch.I) that
\begin{eqnarray*}
 |K_{1\otimes2}| &\le& |K_{1}\otimes I_{2}| + |I_{1}\otimes K_{2}| \\
                 &\le& \sqrt{r_{2}}|K_{1}| + \sqrt{r_{1}}|K_{2}|\,.
\end{eqnarray*}
Since $K_{1}$ and $K_{2}$ are $L^{1}$-bounded, by integrating this inequality over $U_{1}\cap U_{2}$ it follows that 
$K_{1\otimes 2}$ is also $L^{1}$-bounded. Similarly, for the Chern curvature $R_{1\otimes2}$ we obtain
\begin{eqnarray*}
 |R_{1\otimes2}|^{2} &\le& |R_{1}\otimes I_{2}|^{2} + |I_{1}\otimes R_{2}|^{2} + 2|R_{1}\otimes I_{2}||I_{1}\otimes R_{2}| \\
                     &\le& r_{2}|R_{1}|^{2} + r_{1}|R_{2}|^{2} + 2\sqrt{r_{1}r_{2}}|R_{1}||R_{2}|\,.
\end{eqnarray*}
Now, since $R_{1}$ and $R_{2}$ are square-integrable and the product $|R_{1}||R_{2}|$ is $L^{1}$-bounded, the square-integrability of 
$R_{1\otimes2}$ follows integrating the above inequility over $U_{1}\cap U_{2}$. \\

On the other hand, if $R_{1\oplus 2}$ (resp. $K_{1\oplus 2}$) denotes the Chern curvature (resp. the associated mean curvature) of 
$h_{1}\oplus h_{2}$, we obtain
\begin{equation}
 |K_{1\oplus2}| \le |K_{1}| + |K_{2}|\,,  \nonumber
\end{equation}
\begin{equation}
 |R_{1\oplus 2}|^{2} \le |R_{1}|^{2} + |R_{2}|^{2} + 2|R_{1}||R_{2}|\,,  \nonumber
\end{equation}
and the result follows from these inequalities in a similar way to the tensor product case. \;\;Q.E.D. \\

Following \cite{Biswas-Schumacher}, we say that an admissible structure $h$ on a torsion-free Higgs sheaf $\mathfrak E$ is an 
{\it Hermitian-Yang-Mills structure}, if on the open set $U\subset X$ where the metric is defined we have 
\begin{equation}
{\mathcal K}_{h} = K_{h} + i\Lambda [\phi,\bar\phi_{h}] = c \cdot I\,
\end{equation}
where $c$ is the constant used in the Donaldson functional (see \cite{Cardona}). It is important to note here that the admissibility of a 
metric on a Higgs sheaf depends only on conditions imposed on the Chern curvature. In contrast, the notion of Hermitian-Yang-Mills structure 
does depend on the Higgs field conditions imposed on the Hitchin-Simpson curvature. \\

We say that a torsion-free Higgs sheaf $\mathfrak E$ has an {\it approximate Hermitian-Yang-Mills structure}, if for all $\epsilon>0$, there 
exists an admissible structure $h_{\epsilon}$ such that 
\begin{equation}
 \sup_{U_{\epsilon}}|{\cal K}_{\epsilon}-cI|<\epsilon\,.  \label{Cond. approx. HYM structures for Higgs sheaves}
\end{equation}
Here $U_{\epsilon}$ is the open set where $h_{\epsilon}$ is defined and ${\cal K}_{\epsilon}$ is the Hitchin-Simpson mean curvature of 
$h_{\epsilon}$. \\

On the other hand, in an analog way to the classical case, the Hitchin-Simpson mean curvature of the tensor product 
satisfies
\begin{equation}
 |{\cal K}_{1\otimes2}-cI| \le \sqrt{r_{2}}|{\cal K}_{1}-c_{1}I_{1}| + \sqrt{r_{1}}|{\cal K}_{2}-c_{2}I_{2}|\,, \label{ineq. HS curv}
\end{equation}
where $I=I_{1}\otimes I_{2}$ and $c=c_{1}+c_{2}$. Evenmore, if $c=c_{1}=c_{2}$, the Hitchin-Simpson mean curvature of the Whitney sum 
satisfies
\begin{equation}
 |{\cal K}_{1\oplus2}-cI| \le |{\cal K}_{1}-c_{1}I_{1}| + |{\cal K}_{2}-c_{2}I_{2}|\,, \label{ineq. HS curv 2} 
\end{equation} 
where this time $I=I_{1}\oplus I_{2}$. From the inequalities (\ref{ineq. HS curv}), (\ref{ineq. HS curv 2}) and Proposition 
\ref{tensor product of admissible h's} we conclude the following, which is indeed an extension of Proposition 3.1 in \cite{Cardona} 
to torsion-free Higgs sheaves
\begin{pro}
Let ${\mathfrak E}_{1}$ and ${\mathfrak E}_{2}$ be two torsion-free Higgs sheaves over a compact K\"ahler manifold $X$. If they admit approximate 
Hermitian-Yang-Mills structures, so does their tensor product ${\mathfrak E}_{1}\otimes{\mathfrak E}_{2}$. Furthermore, if 
$\mu({\mathfrak E}_{1})=\mu({\mathfrak E}_{2})$, so does their Whitney sum ${\mathfrak E}_{1}\oplus{\mathfrak E}_{2}$. 
\end{pro}

Let $\mathfrak E=(E,\phi)$ be a torsion-free Higgs sheaf over $X$. By Biswas and Schumacher \cite{Biswas-Schumacher} 
(see also \cite{Bando-Siu} or \cite{Hironaka}), there exists a finite sequence of blowups with smooth centers 
\begin{equation}
 \pi_{j}:X_{j}\longrightarrow X_{j-1}\,,  \nonumber
\end{equation}
with $j=1,...,k$ and $X_{0}=X$, such that the pullback of the sheaf $E^{*}$ to $X_{k}$ modulo torsion is locally free and 
$\pi_{1}\cdots\pi_{k}$ outside $S$ is a biholomorphism. In other words, setting $\tilde X=X_{k}$ and
\begin{equation}
 \pi=\pi_{1}\cdots\pi_{k}:\tilde X \longrightarrow X\,,   \label{reg. morphism} 
\end{equation}
and denoting by $T$ the torsion part of $\pi^{*}E^{*}$, then $\pi^{*}E^{*}/T$ is a holomorphic bundle over $\tilde X$ and 
$\pi$ restricted to $\tilde X\backslash\pi^{-1}(S)$ is a biholomorphism. \\

Let $\tilde E$ be the dual of the bundle $\pi^{*}E^{*}/T$. Clearly, the morphism $\phi$ defines a Higgs field 
$\psi=\pi^{*}\phi^{*}$ on $\pi^{*}E^{*}$ and since $\psi(T)\subset T\otimes\Omega^{1}_{\tilde X}$, the morphism $\psi$ is well defined 
on the quotient $\pi^{*}E^{*}/T$ and we have a morphism
\begin{equation}
 \psi^{*}:{\tilde E}\longrightarrow{\tilde E}\otimes\Omega^{1}_{\tilde X}\,.  \nonumber
\end{equation}

From the above analysis we conclude that $\tilde{\mathfrak E}=(\tilde E,\psi^{*})$ is a Higgs bundle over $\tilde X$. We say that 
$\tilde{\mathfrak E}$ is a {\it regularization} of the Higgs sheaf $\mathfrak E$ and that the map $\pi$, defined by (\ref{reg. morphism}), is a 
morphism regularizing $\mathfrak E$. \\

If $\omega$ is a K\"ahler metric on $X$, its pullback $\pi^{*}\omega$ is degenerate along the exceptional divisor $\pi^{-1}(S)$ and hence it 
is not a K\"ahler metric on $\tilde X$. From \cite{Biswas-Schumacher} (see also \cite{Buchdahl-0}), we know there exists a K\"ahler metric closely 
related to the form $\pi^{*}\omega$. This metric can be defined as follows: Let $\eta$ be an arbitrary K\"ahler metric on $\tilde X$ and 
$0\le\epsilon\le 1$, then we define
\begin{equation}
 \omega_{\epsilon}=\pi^{*}\omega + \epsilon\eta\,.   \label{metric omega-epsilon}
\end{equation}
Clearly, this is a K\"ahler metric for each $\epsilon>0$. Such a metric can be used to prove some simple properties involving admissible metrics. In particular, 
we have the following result
\begin{pro}\label{reg. metrics vs admiss. metrics}
 Let $\mathfrak E$ be a torsion-free Higgs sheaf over a compact K\"ahler manifold $X$ and $\tilde{\mathfrak E}$ a regularization of it.
Then, any Hermitian metric on $\tilde{\mathfrak E}$ induces an admissible metric on $\mathfrak E$. 
\end{pro}
 \noindent {\it Proof:} 
Let $\tilde h$ be an Hermitian metric on $\tilde{\mathfrak E}$ and denote by $S$ the singularity set of $\mathfrak E$. Let $\pi$ be the morphism 
regularizing $\mathfrak E$. The Hermitian metric $\tilde h$ induces an Hermitian metric $h$ on ${\mathfrak E}|_{X\backslash S}$. \\

Let $K=i\Lambda R$ be the (classical) mean curvature of ${\mathfrak E}|_{X\backslash S}$ associated with the metric $h$. 
Since $\tilde h$ is defined on all $\tilde X$, the pullback of $R$, denoted by $\tilde R$, extends to all $\tilde X$ as the curvature of the 
Hermitian metric $\tilde h$ and hence, on each point of $\tilde X$, we have 
\begin{equation}
 in{\tilde R}\wedge\omega_{\epsilon}^{n-1} =  i\Lambda_{\epsilon}{\tilde R}\,\omega_{\epsilon}^{n} = {\tilde K}_{\epsilon}\,\omega_{\epsilon}^{n} 
\label{formula K and R} 
\end{equation}
where $i\Lambda_{\epsilon}$ denotes this time the adjoint of the multiplication by $\omega_{\epsilon}$ and ${\tilde K}_{\epsilon}$ represents the 
corresponding mean curvature of $\tilde{\mathfrak E}$. Now, since $\tilde X$ is compact, for some positive constant $C$ 
(see \cite{Bando-Siu}, Lemma 5) we have 
\begin{equation}
 i{\tilde R}\le iC\omega_{\epsilon}I  \label{estimative R - omega}
\end{equation}
where $I$ here is the identity endomorphism of $\tilde E$. Hence, by applying $i\Lambda_{\epsilon}$ to (\ref{estimative R - omega}) we have that 
$C i\Lambda_{\epsilon}\,\omega_{\epsilon}I - {\tilde K}_{\epsilon}$ must be a semipositive definite endomorphism of $\tilde E$ and we get
\begin{eqnarray*}
 |{\tilde K}_{\epsilon}|\,\omega_{\epsilon}^{n} &\le& |{\tilde K}_{\epsilon} - Ci\Lambda_{\epsilon}\omega_{\epsilon}I|\,\omega_{\epsilon}^{n} +  
                                           |Ci\Lambda_{\epsilon}\omega_{\epsilon}I|\,\omega_{\epsilon}^{n}\\
                                      &\le& {\rm tr}\,[i\Lambda_{\epsilon}(C\omega_{\epsilon}I - {\tilde R})]\,\omega_{\epsilon}^{n} + 
                                           {\rm tr}\,(Ci\Lambda_{\epsilon}\omega_{\epsilon}I)\,\omega_{\epsilon}^{n}\\
                                      &\le& in\,{\rm tr}\,(2C\omega_{\epsilon}I - {\tilde R})\,\omega_{\epsilon}^{n-1}\\
                                      &\le& in\,{\rm tr}\,(2C\omega_{1}I - {\tilde R})\,\omega_{1}^{n-1}\,.
\end{eqnarray*}                            
From this we conclude that ${\tilde K}_{\epsilon}$ is uniformly integrable with respect to $0<\epsilon\le 1$  and hence, taking the limit 
$\epsilon\rightarrow 0$, it follows that $K$ is $L^{1}$-bounded. \\

On the other hand, from the theory of holomorphic bundles (see in particular \cite{Bando-Siu}, Lemma 6) we have the following identity 
\begin{equation}
\left[2\,c_{2}(\tilde{E}) - c_{1}(\tilde{E})^{2}\right]\cup[\omega_{\epsilon}]^{n-2}  =  
                \frac{1}{4\pi^{2}n(n-1)}\int_{\tilde X}\left[ |{\tilde R}|^{2} - |{\tilde K}_{\epsilon}|^{2} \right]\omega_{\epsilon}^{n}\,,
\end{equation}
and therefore, taking the limit $\epsilon\rightarrow 0$, it follows that $R$ is also square-integrable.  \;\;Q.E.D. \\

Notice that since the codimension of the singularity set $S$ of a torsion-free Higgs sheaf $\mathfrak E$ is greater or equal than two, we can 
see ${\mathfrak E}$ as a Higgs bundle over the non-compact manifold ${X\backslash S}$. Thus, in studying torsion-free Higgs sheaves we are 
considering implicitly Higgs bundles over non-compact K\"ahler manifolds. Evenmore, from Section 7 in \cite{Cardona} we know that 
$X\backslash S$ satisfies all assumptions that Simpson \cite{Simpson} imposes on the base manifold and that the Hitchin-Simpson mean 
curvature ${\cal K}$ is $L^{1}$-bounded on $X\backslash S$. Clearly, these properties extend naturally to the open sets $U$ where the 
admissible structures are defined.

\section{Higgs sheaves and admissible structures}

As we said before, Biswas and Schumacher proved a Hitchin-Kobayashi correspondence for polystable Higgs sheaves. As an immediate consequence of this
we know that any restriction of a stable Higgs sheaf $\mathfrak E$ to $X\backslash S$ is $\omega$-polystable. In fact, using the modified 
definition of admissibility, this also holds for restrictions to certain open subsets of $X$. To be precise we have the following result
\begin{pro}\label{s sheaf ps restriction}
 Let $\mathfrak E$ be a torsion-free Higgs sheaf over a compact K\"ahler manifold $X$ with K\"ahler form $\omega$ and denote by $S$ its 
singularity set. Let $U\subset X$ be an open set such that $X\backslash U$ is a closed analytic subset of codimension at least two containing 
$S$. Then ${\mathfrak E}|_{U}$ is $\omega$-polystable if $\mathfrak E$ is $\omega$-polystable. 
\end{pro}
 \noindent {\it Proof:} 
Let $\mathfrak E$ be a torsion-free sheaf over $X$ and assume first that it is $\omega$-stable. Then, by Theorem 3.1 in \cite{Biswas-Schumacher}, 
there exists an Hermitian-Yang-Mills structure $h$ on it. Let $U$ be an open subset of $X$ such that $X\backslash U$ is a closed analytic subset 
with codimension greater or equal than two and suppose that $S\subset X\backslash U$. Thus, $h$ is in particular an Hermitian-Yang-Mills metric 
on the Higgs bundle ${\mathfrak E}|_{U}$ and hence, from Proposition 3.3 in \cite{Simpson}, it must be $\omega$-polystable.  \\
 
Assume now that $\mathfrak E$ is $\omega$-polystable. Then, it can be decomposed as a direct sum of $\omega$-stable Higgs sheaves with the same 
slope as $\mathfrak E$. From the first part of the proof, we know that each restriction of these stable Higgs sheaves to $U$ must be 
$\omega$-polystable and hence the result follows.  \;\;Q.E.D. \\

\begin{lem}\label{Lemma polystability}
 Let ${\mathfrak E}_{1}$ and ${\mathfrak E}_{2}$ be two torsion-free Higgs sheaves over a compact K\"ahler manifold $X$ with K\"ahler form $\omega$.
If both are $\omega$-polystable, then ${\mathfrak E}_{1}\otimes{\mathfrak E}_{2}$ modulo torsion is also $\omega$-polystable. 
\end{lem}

\noindent {\it Proof:} Let ${\mathfrak E}_{1}$ and ${\mathfrak E}_{2}$ be $\omega$-polystable. Then, from Corollary 3.5 in
\cite{Biswas-Schumacher} we know there exist Hermitian-Yang-Mills structures $h_{1}$ and $h_{2}$. Now, by Proposition 
\ref{tensor product of admissible h's}, it follows that $h=h_{1}\otimes h_{2}$ is an admissible structure on 
${\mathfrak E}_{1}\otimes{\mathfrak E}_{2}$. Clearly, it is an Hermitian-Yang-Mills structure and hence, using again 
the same Corollary, such a tensor product (modulo torsion) must be $\omega$-polystable.  \;\;Q.E.D. \\

As a consequence of the above Lemma, Biswas and Schumacher \cite{Biswas-Schumacher} proved that the tensor product of two semistable sheaves is again semistable. Here 
we present a different proof. Notice first that from Lemma \ref{Lemma polystability} we have the following result

\begin{lem}\label{2do Imte Lemma polystability}
 Let ${\mathfrak E}_{1}$ and ${\mathfrak E}_{2}$ be two torsion-free Higgs sheaves over a compact K\"ahler manifold $X$ with K\"ahler form 
$\omega$. If ${\mathfrak E}_{1}$ is $\omega$-semistable and ${\mathfrak E}_{2}$ is $\omega$-polystable, then ${\mathfrak E}_{1}\otimes{\mathfrak E}_{2}$ 
modulo torsion is $\omega$-semistable. 
\end{lem}

\noindent {\it Proof:} Assume that ${\mathfrak E}_{2}$ is $\omega$-polystable an ${\mathfrak E}_{1}$ is $\omega$-semistable. Following 
Simpson \cite{Simpson 3} (see also \cite{Biswas-Schumacher}), there exists a filtration of ${\mathfrak E}_{1}$ by Higgs subsheaves 
\begin{equation}
 0={\mathfrak F}_{0}\subset {\mathfrak F}_{1}\subset\cdots\subset {\mathfrak F}_{k}={\mathfrak E}_{1}\,, \label{HN filtration}
\end{equation}
in which the quotients ${\mathfrak F}_{j}\slash {\mathfrak F}_{j-1}$ for $j=1,...,k$ are $\omega$-polystable and they have all the 
same slope as ${\mathfrak E}_{1}$. Now, let $U\subset X$ be the open set in which all terms of the filtration (\ref{HN filtration}), 
all quotients ${\mathfrak F}_{j}\slash {\mathfrak F}_{j-1}$, and also ${\mathfrak E}_{2}$ are locally free. Then $X\backslash U$ is a 
closed analytic subset of codimension greater or equal than two, and on $U$ we have the sequence 
\begin{equation}
 \xymatrix{
0 \ar[r]  &  {\mathfrak F}_{1} \ar[r]  &   {\mathfrak F}_{2} \ar[r]  &   {\mathfrak F}_{2}/{\mathfrak F}_{1}  \ar[r]  &    0  
} \label{HN filtration 1}
\end{equation}
as a sequence of locally free Higgs sheaves. Then, tensoring the above sequence by ${\mathfrak E}_{2}$ we obtain the sequence   
\begin{equation}
 \xymatrix{
0 \ar[r]  &  {\mathfrak F}_{1}\otimes{\mathfrak E}_{2} \ar[r]  &  {\mathfrak F}_{2}\otimes{\mathfrak E}_{2} \ar[r]  &  ({\mathfrak F}_{2}/{\mathfrak F}_{1})\otimes{\mathfrak E}_{2} \ar[r]  &    0  
} \label{HN filtration 2}
\end{equation}
which is again an exact sequence of locally free Higgs sheaves over $U$. Since ${\mathfrak F}_{1}$ and ${\mathfrak F}_{2}/{\mathfrak F}_{1}$ are 
both $\omega$-polystable, and also ${\mathfrak E}_{2}$ is $\omega$-polystable by hypothesis, we have by Proposition \ref{s sheaf ps restriction} that they are 
all $\omega$-polystable over $U$. Therefore, it follows from Lemma \ref{Lemma polystability} that ${\mathfrak F}_{1}\otimes{\mathfrak E}_{2}$ and 
$({\mathfrak F}_{2}/{\mathfrak F}_{1})\otimes{\mathfrak E}_{2}$ are both $\omega$-polystable with equal slopes (in particular they are 
$\omega$-semistable). Therefore, from this and Corollary \ref{extension of ss sheaves}, we obtain the semistability of the Higgs sheaf 
${\mathfrak F}_{2}\otimes{\mathfrak E}_{2}$ over the open set $U$. \\

Now, we consider the exact sequence
\begin{equation}
 \xymatrix{
0 \ar[r]  &  {\mathfrak F}_{2} \ar[r]  &   {\mathfrak F}_{3} \ar[r]  &   {\mathfrak F}_{3}/{\mathfrak F}_{2}  \ar[r]  &    0\,.  
} 
\end{equation}
Since over $U$ this is an exact sequence of locally free Higgs sheaves, tensoring again by ${\mathfrak E}_{2}$ we obtain over $U$ the following 
exact sequence of locally free Higgs sheaves:
\begin{equation}
 \xymatrix{
0 \ar[r]  &  {\mathfrak F}_{2}\otimes{\mathfrak E}_{2} \ar[r]  &   {\mathfrak F}_{3}\otimes{\mathfrak E}_{2} \ar[r]  &  ({\mathfrak F}_{3}/{\mathfrak F}_{2})\otimes{\mathfrak E}_{2} \ar[r]  &    0\,.  
} 
\end{equation}
Using again Lemma \ref{Lemma polystability} we have that $({\mathfrak F}_{3}/{\mathfrak F}_{2})\otimes{\mathfrak E}_{2}$ is $\omega$-polystable, in 
particular it is $\omega$-semistable and since ${\mathfrak F}_{2}\otimes{\mathfrak E}_{2}$ is also $\omega$-semistable, we obtain (again by Corollary 
\ref{extension of ss sheaves}) that ${\mathfrak F}_{3}\otimes{\mathfrak E}_{2}$ is $\omega$-semistable. Continuing this process we get at 
the end that ${\mathfrak E}_{1}\otimes{\mathfrak E}_{2}$ is $\omega$-semistable. Since all of this holds over $U$, whose complement has 
codimension greater or equal than two, it can be extended on all $X$ and ${\mathfrak E}_{1}\otimes{\mathfrak E}_{2}$ is $\omega$-semistable 
on $X$ as well. \;\;Q.E.D. \\

\begin{thm}\label{Th. tensor ss}
  Let ${\mathfrak E}_{1}$ and ${\mathfrak E}_{2}$ be two torsion-free Higgs sheaves over a compact K\"ahler manifold $X$ with K\"ahler form $\omega$.
If both are $\omega$-semistable, then ${\mathfrak E}_{1}\otimes{\mathfrak E}_{2}$ modulo torsion is also $\omega$-semistable. 
\end{thm}

\noindent {\it Proof:} The Higgs sheaf ${\mathfrak E}_{1}$ has a filtration by Higgs subsheaves as in (\ref{HN filtration}), with 
$\omega$-polystable quotients with the same slope as ${\mathfrak E}_{1}$. Now, let $U\subset X$ be an open subset such that all terms of the 
filtration, all quotients and also ${\mathfrak E}_{2}$ are locally free. Then, we have the exact sequences (\ref{HN filtration 1}) and 
(\ref{HN filtration 2}) and since ${\mathfrak E}_{2}$ is $\omega$-semistable, the result follows by applying Lemma \ref{2do Imte Lemma polystability}. 
 \;\;Q.E.D. \\


\begin{thebibliography}{99}

\bibitem{Atiyah-Bott} {\sc Atiyah, M., Bott, R.,} {\it The Yang-Mills equations over Riemann surfaces,} Phil. Trans. Roy. Soc. London A., {\bf 308}  (1983), pp. 532-615.

\bibitem{Bando-Siu} {\sc Bando, S., Siu,  Y.-T.,} \emph{Stable sheaves and Einstein-Hermitian metrics}, Geometry and analysis on complex manifolds, (World Scientific 
Publishing, River Edge, NJ, 1994), pp. 39-50.

\bibitem{Biswas-Schumacher} {\sc Biswas, I., Schumacher, G.,} {\it Yang-Mills equations for stable Higgs sheaves}, International Journal of Mathematics, Vol. 20, No 5 (2009), 
pp. 541-556.

\bibitem{Bruzzo-Granha} {\sc Bruzzo, U., Gra\~na Otero, B.,} {\it Metrics on semistable and numerically effective Higgs bundles}, J. reine ang. Math.,  {\bf 612} (2007), pp. 59-79.

\bibitem{Buchdahl-0} {\sc Buchdahl, N. M.,} {\it Hermitian-Einstein connections and stable vector bundles over compact complex surfaces}, Math. Ann. {\bf 280} (1988), pp. 625-648.

\bibitem{Cardona} {\sc Cardona, S. A. H.,} \emph{Approximate Hermitian-Yang-Mills structures and semistability for Higgs bundles. I: generalities 
and the one-dimensional case}, Ann. Glob. Anal. Geom., Vol. 42, Number 3 (2012), pp. 349-370 (DOI 10.1007/s10455-012-9316-2).

\bibitem{Donaldson-1} {\sc Donaldson, S. K.,} {\it A new proof of a theorem of Narasimhan and Seshadri,} J. Diff. Geom. Soc., {\bf 18}  (1983), pp. 269-278.

\bibitem{Griffiths-Harris} {\sc Griffiths, P., Harris, J.,} \emph{Principles of algebraic geometry}, John Wiley and Sons. (1978).

\bibitem{Hironaka} {\sc Hironaka, H.,} {\it Flattening theorem in complex-analytic geometry}, American Journal of Mathematics, Vol. 97, No. 2 (Summer, 1975), 
pp. 503-547.

\bibitem{Hitchin} {\sc Hitchin, N. J.,} {\it The self-duality equations on a Riemann surface}, Proc. London. Math., {\bf 55} (1987), pp. 59-126.

\bibitem{Kobayashi} {\sc Kobayashi, S.,} \emph{Differential geometry of complex vector bundles}, Publications of the Mathematical Society of Japan, Vol. 15 
(Iwanami Shoten Publishers and Princeton University Press, 1987). 

\bibitem{Lubke} {\sc L\"ubke, M., Teleman, A.,} \emph{The Kobayashi-Hitchin correspondence}, World Scientific Publishing Co. Pte. Ltd., (1995).

\bibitem{Narasimhan-Seshadri} {\sc Narasimhan, M., Seshadri, C.,} {\it Stable and unitary bundles on a compact Riemann surface,} Math. Ann., {\bf 82}  (1965), pp. 540-564.

\bibitem{Simpson} {\sc Simpson, C. T.,} \emph{Constructing variations of Hodge structure using Yang-Mills theory and applications to uniformization,} 
J. Amer. Math. Soc. {\bf 1} (1988), pp. 867-918. 

\bibitem{Simpson 2} {\sc Simpson, C. T.,} {\it Higgs bundles and local systems,} Publ. Math. I.H.E.S., {\bf 75} (1992), pp. 5-92.

\bibitem{Simpson 3} {\sc Simpson, C. T.,} {\it Moduli of representations of the fundamental group of a smooth projective variety I,} Publ. Math. I.H.E.S. {\bf 79} (1994), pp. 47-129.

\bibitem{Siu} {\sc Siu, Y.-T.,} \emph{Lectures on Hermitian-Einstein metrics for stable bundles and K\"ahler-Einstein metrics}, Birkh\"auser, Basel-Boston, (1987).

\bibitem{Uhlenbeck-Yau} {\sc Uhlenbeck, K., Yau, S.-T.,} {\it On the existence of Hermitian-Yang-Mills connections in stable vector bundles,} Commun. Pure Appl. Math., 
{\bf 39}  (1986), pp. S257-S293.

\bibitem{Jiayu-Zhang} {\sc Li, J., Zhang, X.,} \emph{Existence of approximate Hermitian-Einstein structures on 
semistable Higgs bundles,} arXiv:1206.6676v1, 2012.


\end{thebibliography}
\end{document}